\documentclass[11pt,leqno]{article}
\usepackage{amsmath,amssymb,amsfonts}

\newtheorem{theorem}{Theorem}[section]
\newtheorem{proposition}[theorem]{Proposition}
\newtheorem{corollary}[theorem]{Corollary}
\newtheorem{lemma}[theorem]{Lemma}
\newtheorem{definition}[theorem]{Definition}

\newtheorem{question}[theorem]{Question}

\def\IP{\hbox{\bf IP}}
\def\DMO{\hbox{\bf DMO}}

\newenvironment{proof}{\mbox{\bf Proof.}}{\mbox{$\blacksquare$}\bigskip}

\begin{document}

\begin{center}

{\LARGE\bf The Skolem-Bang Theorems

in Ordered Fields with an $\IP$}

\vspace{4mm}

{\large{\bf Seyed Masih Ayat}}\hspace{5mm} 

 Institute for Studies in Theoretical Physics and Mathematics
(IPM), Tehran, Iran

and

\noindent Department of Mathematics, Tarbiat Modarres University,
Tehran, Iran

smasih@ipm.ir

\end{center}

\begin{abstract}
This paper is concerned with the extent to which the Skolem-Bang
theorems in Diophantine approximations generalise from the
standard setting of $\langle\mathbb{R},\mathbb{Z}\rangle$ to
structures of the form $\langle F,I\rangle$, where $F$ is an
ordered field and $I$ is an integer part of $F$. We show that some
of these theorems are hold unconditionally in general case
(ordered fields with an integer part). The remainder results are
based on Dirichlet's and Kronecker's theorems. Finally we extend
Dirichlet's theorem to ordered fields with $IE_1$ integer part.

\end{abstract}

{\large\bf{Subj-Class: LO}}


{\large\bf{keywords:}} Nonstandard Diophantine Approximations, The
Skolem-Bang Integer Part Theorems, Dense Mod 1 Subsets, Weak
Fragments of Arithmetic.

\section{Introduction }
Let $\alpha\geqslant 1$ be a real number. The notion of
$\mathbb{N}_\alpha$ was introduced by Skolem and Bang as the
sequence $\{\lfloor n\alpha\rfloor|\ n\in\mathbb{N}\ \}$ of
positive integers, where $\lfloor x\rfloor$ is the integer part of
$x$. The following facts are studied in Skolem-Bang Theorems
\cite{Sk,Ba}:

\noindent {\bf 1}. $\mathbb{N}_\alpha\cap\mathbb{N}_\beta=\{0\}$;
{\bf 2}. $\mathbb{N}_\alpha\cup\mathbb{N}_\beta=\mathbb{N}$; {\bf
3}. $\mathbb{N}_\alpha\subseteq\mathbb{N}_\beta$.

\noindent These theorems are also reported in \cite{N}:

\noindent {\bf Fact A}. Let $\alpha$ and $\beta$ be positive real
numbers. Then $\mathbb{N}_\alpha\cup\mathbb{N}_\beta=\{0\}$ and
$\mathbb{N}_\alpha\cap\mathbb{N}_\beta=\{0\}$ if and only if
$\alpha,\beta$ are irrational numbers and
$\alpha^{-1}+\beta^{-1}=1$.

\noindent {\bf Fact B}. Let $\alpha$ and $\beta$ be positive real
numbers. If $1,\alpha^{-1},\beta^{-1}$ are linearly independent over
the field of rational numbers, then $\mathbb{N}_\alpha$ and
$\mathbb{N}_\beta$ have infinitely many common elements.

\noindent {\bf Fact C}. Let $\alpha$ and $\beta$ be positive real
numbers such that $a\alpha^{-1}+b\beta^{-1}=c$ for some integers
$a,b,c$, with $ab<0$ and $c\neq0$. Then $\mathbb{N}_\alpha$ and
$\mathbb{N}_\beta$ have infinitely many common elements.

\noindent {\bf Fact D}. Let $\alpha$ and $\beta$ be positive real
numbers such that $a\alpha^{-1}+b\beta^{-1}=c$ for some positive
integers $a,b,c$, with $(a,b,c)=1$ and $c>1$. Then
$\mathbb{N}_\alpha$ and $\mathbb{N}_\beta$ have infinitely many
common elements.

\noindent {\bf Fact E}. Let $\alpha$ and $\beta$ be positive real
numbers. The sets $\mathbb{N}_{\alpha}$ and $\mathbb{N}_{\beta}$
are disjoint if and only if $\alpha$ and $\beta$ are irrational
numbers and there exist positive integers $a$ and $b$ such that
$a\alpha^{-1}+b\beta^{-1}=1$.

\noindent Further more if $\mathbb{N}_\alpha$ and
$\mathbb{N}_\beta$ have one common element, they have infinitely
many ones.

\noindent {\bf Fact F}. Let $\alpha$ and $\beta$ be positive
irrational numbers. Then
$\mathbb{N}_\alpha\supseteq\mathbb{N}_\beta$ if and only if there
exist positive integers $a$ and $b$ such that
$a(1-\alpha^{-1})+b\beta^{-1}=1$.

\noindent The rational version of Fact F is the following:

\noindent {\bf Fact F'}. Let $\sigma$ and $\rho$ be positive
rational numbers. Then $\mathbb{N}_\sigma\supseteq\mathbb{N}_\rho$
if and only if there exist positive integers $a$ and $b$ such that
$a(1-\sigma^{-1})+b\rho^{-1}=1$.

There are also some other Skolem-Bang results which are either
trivial or obtained from the above ones. All these results are
based on two important theorems in the theory of Diophantine
Approximations: {\it Dirichlet's Theorem} and {\it Kronecker's
Theorem}.

\noindent {\bf Dirichlet's Theorem}. Let $\theta$ be a positive
irrational number. There are infinitely many rational numbers
$\frac{a}{b}$, where $a$ and $b$ are positive integers, such that
$$|\theta-\frac{a}{b}|<\frac{1}{b^2}\ .$$

\noindent An immediate conclusion of Dirichlet's Theorem is that
the set $\{n\theta-\lfloor n\theta\rfloor|\ n\in\mathbb{N}\}$ is a
dense subset of $[0,1)$. A more interesting corollary is

\noindent {\bf Separability property}. Let $\alpha,\beta>1$ be
real numbers. Then $\alpha\neq\beta$ if and only if
$\mathbb{N}_\alpha\neq\mathbb{N}_\beta$.

\noindent {\bf Kronecker's Theorem }. Let $\alpha$ and $\beta$ be
positive irrationals such that $1,\alpha,\beta$ are linearly
independent over the field of rational numbers, then the points
whose coordinates are the fractional parts of multiples of
$\alpha$ and $\beta$, i.e. $(n\alpha-\lfloor
n\alpha\rfloor,n\beta-\lfloor n\beta\rfloor),\ n=1,2,3,\dots$, are
dense in the {\it unite square}.

Note that Fact B and Kronecker's Theorem are equivalent.
Dirichlet's Theorem and Kronecker's Theorem are based on Pigeon
Hole Principle ({\bf PHP}) and Box Principle (the two dimensional
version of Pigeon Hole Principle). However, in non-Archimedean
cases, PHP and Box Principle do not hold. Extending the notion of
the {\it separability property} to non-Archimedean structures
$\langle F,I\rangle$ is a useful tool to generalize Skolem-Bang
Theorems to these structures. If $F$ is an ordered field and $I$
an integer part for $F$, we call $\langle F,I\rangle$ {\it
separable} if it satisfies the separability property. Mojtaba
Moniri has conjectured that ``{\it any arbitrary structure
$\langle F,I\rangle$ is separable}". In Section 2, we prove some
weak versions of the separability property, i.e. we prove it for
the cases that:
\begin{enumerate}
\item $\alpha$ and $\beta$ are irrationals;
\item $\alpha$ and $\beta$ are rationals;
\item $\alpha,\beta\geqslant 2$;
\item $\rho$ is a rational number and $\alpha$ is an irrational
number such that $1<\rho<\alpha<2$.
\end{enumerate}

In Section 3, we will show that any $\langle F,I\rangle$
satisfying Dirichlet's Theorem is separable and Fact A hold in
separable $\langle F,I\rangle$. Also if $\langle F,I\rangle$
satisfies Dirichlet's Theorem and $I$ is a B\'ezout domain then
Facts C, D and one direction of Facts E, F hold. We also show that
Fact F' holds for a structure $\langle F,I\rangle$ in which $I$ is
a B\'ezout domain.

The main tool of Section 4 is {\it Farey series} which is studied
in Hardy and Wright's excellent book \cite{HW}. Using {\it weak}
versions of {\bf PHP}, we can prove some special forms of
Dirichlet's Theorem in {\it weak} fragments of Arithmetic. In
\cite[Theorem 3.1]{D}, P. D'Aquino proved a weak version of
Dirichlet's Theorem in $I\Delta_0+\Omega_1$, where $\Omega_1$ is
$$\forall x \exists y (x^{log(x)}=y)$$
and by $log(x)$ we mean the integer part of $\log_{2}(x)$, (for
more details, see subsection 4.1). Using Farey series, we prove
Dirichlet's Theorem in $\langle F,I\rangle$ in which $I$ is an
$IE_1$-model. Since Wilmers proved that $IE_1\vDash B\acute{e}z$,
Facts C, D and one direction of Facts E, F mentioned in Section 3
hold in $IE_1$, \cite{Wi}. In \cite{Se}, B. Segre provided an
asymmetric Diophantine approximations theorem for irrational
numbers. We prove this theorem by a similar method based on Farey
series which represented in proof theorem 1.7 in \cite{N}. This
theorem has interesting corollaries such as Hurwitz asymmetric
Theorem which will be denoted in Section 4. Their proofs are
similar to the real case.

It is not known whether Kronecker's Theorem holds for
$IE_1$-models. If so, all the Skolem-Bang theorems hold in any
$\langle F,I\rangle$ in which $I$ is an $IE_1$-model. In fact the
remainder direction of Facts E, F are based on Kronecker's
theorem.

\subsection{The Preliminaries}
Let $L=\{+,-,\cdot,0,1,\leqslant\}$ be the language of ordered
rings. We deal with the following sets of axioms in L:

${\bf DOR}$: discretely ordered rings i.e., axioms for ordered
rings together with $\forall x \neg(0<x<1)$.

${\bf ZR}$: discretely ordered $\mathbb{Z}$-rings i.e., ${\bf
DOR}$ together with the condition that for every
$n\in\mathbb{N}^{\geqslant 2}$, we have $(\forall x)(\exists q,r)
(x=nq+r\bigwedge 0\leqslant r<n)$.

${\bf EDR}$: Euclidean division rings i.e., ${\bf DOR}$ extended
with the scheme of axioms that for every $n\in I^{> 0}$, $(\forall
x)(\exists q,r)(x=nq+r\bigwedge 0\leqslant r<n)$.

${\bf IOP}$: Open induction i.e., ${\bf DOR}$ plus the following
scheme for open $L$-formulas $\varphi(x,y)$

$$\forall \vec{x},(\varphi(\vec{x},0) \bigwedge \forall y\geqslant
0, (\varphi(\vec{x},y)\rightarrow\varphi
(\vec{x},y+1))\rightarrow\forall y\geqslant 0,
\varphi(\vec{x},y)).$$

We define the formula class $E_n, U_n,\forall_n,\exists_n$ in the
usual way:
$$E_0=U_0=\{\phi(\bar{x}):\phi\ is\ open\}$$
$$\exists_{n+1}=\{\exists\bar{y}\phi(\bar{x},\bar{y}):\phi\in\forall_n\},\
\forall_{n+1}=\{\forall\bar{y}\phi(\bar{x},\bar{y}):\phi\in\exists_n\}$$
$$E_{n+1}=\{\exists\bar{y}\leqslant t(\bar{x})\phi(\bar{x},\bar{y}):
\phi\in U_n,\ t\ a\ term\ in\ L\}$$
$$U_{n+1}=\{\forall\bar{y}\leqslant t(\bar{x})\phi(\bar{x},\bar{y}):
\phi\in E_n,\ t\ a\ term\ in\ L\}$$
$$\Delta_0=\bigcup_{n\in\mathbb{N}}E_n=\bigcup_{n\in\mathbb{N}}U_n$$

${\bf IE_{1}}$: Bounded existential induction {i.e., ${\bf DOR}$
plus the induction schema for all $E_{1}$-formulas $\varphi$:

\[\forall \vec{x},(\varphi(\vec{x},0)
\bigwedge \forall y\geqslant 0,
(\varphi(\vec{x},y)\rightarrow\varphi
(\vec{x},y+1))\rightarrow\forall y\geqslant 0,
\varphi(\vec{x},y))).\]

\noindent We can define $IE_n,IU_n,I\Delta_0$ similarly.

We say that a subring $I$ of an ordered field $F$ is an integer
part $({\bf IP})$ of $F$ if $I\vDash{\bf DOR}$ and for every $x\in
F$, there is $a\in I$ such that $a\leqslant x<a+1$. We call this
unique element $a$ the integer part of $x$ and write $a=\lfloor
x\rfloor_{I}$. Every real closed field has an $\IP$, \cite{MR}. On
the other hand, there exist ordered fields without any $\IP$, (see
\cite{Bo,K}). One can see that every $\IP$ is an ${\bf EDR}$ and
every ${\bf EDR}$ is an $\IP$ for its fraction field.

We use $\langle F,I\rangle$, for an ordered field $F$ equipped
with an $\IP$ $I$. We set $Q=Frac(I)$, the fraction field of $I$.
We say $I$ is {\bf B\'ezout} if for each $m,n\in I^{\neq 0}$,
there exist $r,s\in I$ such that $rm+sn\geqslant 1$ and $rm+sn |
m,n$. Thus $rm+sn$ is greatest common divisor of $m,n$.}

\section{The Skolem-Bang Integer Part Theorems}

Skolem and Bang theorems (see \cite{Sk, Ba} and \cite{N}), for the
standard case, is based on very special properties of $\mathbb{R}$
and $\mathbb{Z}$, such as {\bf PHP}. In this section, we deal with
these theorems in our an arbitrary $\langle F,I\rangle$.

Fix $\langle F,I\rangle$ and let  $m,k\in I^{\geqslant 0}$, with
$0\leqslant k<m$ and $\alpha\in\ F^{> 0}$. We define an  {\it
arithmetical progression} by parameters $m$ and $k$ as follows,
$mI^{\geqslant 0}+k=\{mt+k\ |\ t\in I^{\geqslant 0}\}$. As the
classic case, let $N_{\alpha}=\{\lfloor n\alpha\rfloor_{I}\ |\
n\in I^{\geqslant 0}\}$ and $\alpha I^{\geqslant 0}=\{ n\alpha\ |\
n\in I^{\geqslant 0}\}$. It is easy to verify $N_\alpha$ when
$0<\alpha<1$.

\begin{lemma} \label{le2} Let $\alpha\in F$. Then
$0<\alpha\leqslant 1$ if and only if $N_{\alpha}=I^{\geqslant 0}$.
\end{lemma}

\begin{proof} ({\bf Only if.}) Suppose $0<\alpha<1$ and pick an
arbitrary $n\in I^{>0}$. Let $k=\lfloor\frac{n}{\alpha}\rfloor$.
We have $k\leqslant\frac{n}{\alpha}<k+1$ and so
$n<(k+1)\alpha\leqslant n+\alpha$ which is less than $n+1$.
Therefore $\lfloor(k+1)\alpha\rfloor=n\in I^{>0}$.

({\bf If.}) Suppose $\alpha>1$ and let
$k=\lfloor\frac{1}{\alpha-1}\rfloor$. Then
$k\leqslant\frac{1}{\alpha-1}<k+1$. So
\[k<k\alpha\leqslant
k+1<k+2<(k+1)\alpha.\]

\noindent We now distinguish two cases. If $k\alpha<k+1$, then
$\lfloor k\alpha\rfloor=k$ and $\lfloor
(k+1)\alpha\rfloor\geqslant k+2$. Therefore $k+1\not\in
I^{\geqslant 0}_{\alpha}$.

If $k\alpha=k+1$, then $\alpha=\frac{k+1}{k}$. So, there is no
$n\in I^{\geqslant 0}$ such that $\lfloor n\alpha\rfloor=k$. The
reason is that $(k-1)\alpha=k+1-\alpha<k<k+1=k\alpha$.
\end{proof}

\noindent In the following theorem we show that when $\alpha>1$ is a
positive rational number, $N_\alpha$ is a union of some arithmetical
progressions. Moreover, if $\alpha,\beta>1$ are two distinct
rational numbers, then $N_\alpha\cap N_\beta$ and
$I^{>0}\setminus(N_\alpha\cup N_\beta)$ are cofinal subsets of
$I^{>0}$.

\begin{theorem} \label{rational} We have
\begin{enumerate}
\item If $\alpha=\frac{p}{q}>1$, $ p,q\in I^{\geqslant 0}$ and $q\neq
0$, then $ N_{\alpha}= \bigcup_{0\leqslant r<q}(pI^{\geqslant
0}+\lfloor\frac{pr}{q}\rfloor)$ and
\[(pI^{\geqslant 0}+(p-1))\subseteq I^{\geqslant 0}\setminus
N_{\alpha}.\]

\item If $\alpha,\beta >1$ are rationals, then $N_{\alpha}\cap
N_{\beta}$ and $ I^{\geqslant 0}\setminus(N_{\alpha}\cup N_{\beta})$
are cofinal subsets of I.

\item If $\alpha_{i}>1$, $i=1,\cdots,n$ are rationals, then
$\bigcap_{i=1}^{n} N_{\alpha_{i}}$ and $ I^{\geqslant
0}\setminus\bigcup_{i=1}^{n} N_{\alpha_{i}}$ are cofinal in I.
\end{enumerate}
\end{theorem}

\begin{proof} 1) Since $I$ is an Euclidean division ring ({\bf EDR}), for each
$n\in I^{\geqslant 0}$ there exist $r,k\in I^{\geqslant 0}$ such
that $0\leqslant r<q$ and $n=kq+r$. Therefore
$n\frac{p}{q}=(kq+r)\frac{p}{q}=kp+\frac{pr}{q}$. Then
$\frac{pr}{q}\leqslant\frac{p(q-1)}{q}<p-1$. So,
$N_{\alpha}=N_{\frac{p}{q}}=\bigcup^{\circ}_{0\leqslant
r<q}(pI^{\geqslant 0}+\lfloor\frac{pr}{q}\rfloor)$. If  $r,s\in
I^{\geqslant 0}$, where $0\leqslant r<s<q$, then we have
$\frac{ps}{q}-\frac{pr}{q}=\frac{p(s-r)}{q}\geqslant
\frac{p}{q}>1$. Therefore
$\lfloor\frac{pr}{q}\rfloor\not\equiv\lfloor\frac{ps}{q}\rfloor(mod\
p)$, since $0\leqslant\frac{pr}{q}<\frac{ps}{q}<p$. So, these
arithmetical progressions are disjoint. The other arithmetical
progressions modulo p appear in $I^{\geqslant 0}\setminus
N_{\alpha}$ as $pI^{\geqslant 0}+(p-1)$.

2) Let $\alpha=\frac{p_1}{q_1},\beta=\frac{p_2}{q_2}$, with
$p_i,q_i\in I^{\geqslant 0},i=1,2$. Then $p_1p_2I^{\geqslant
0}\subseteq N_{\alpha}\cap N_{\beta}$ and
\[(p_1p_2I^{\geqslant 0}+(p_1p_2-1))\cap (N_{\alpha}\cup N_{\beta})=\emptyset.\]

3) Let $\alpha_i=\frac{p_i}{q_i}$, with $p_i,q_i\in I^{\geqslant
0},i=1,2,...,n$. Then we have $(\pi_{i=1}^{n} p_i)I^{\geqslant
0}\subseteq \bigcap_{i=1}^{n} N_{\alpha_i}$ and $[((\pi_{i=1}^{n}
p_i)I^{\geqslant 0}+((\pi_{i=1}^{n} p_i)-1))]\cap
(\bigcup_{i=1}^{n}N_{\alpha_i})=\emptyset$. This completes the
proof.
\end{proof}

We have different situations for $N_\alpha$ with respect to rational
and irrational elements when $\alpha>1$. First we prove a basic
property when $\alpha,\beta\geqslant 2$.

\begin{theorem}\label{T2} Let $\alpha,\beta\in F$ with $\alpha>\beta\geqslant
2$. Then $N_{\beta}\setminus N_{\alpha}\neq\emptyset$.
\end{theorem}

\begin{proof} Let $m=\lfloor\frac{1}{\alpha-\beta}\rfloor$.
If $m=0$, then $\alpha-\beta>1$ and so $\lfloor\beta\rfloor\in
N_{\beta}\setminus N_{\alpha}$. Otherwise, $m>0$. We claim that
$\lfloor (m+1)\beta\rfloor\in N_{\beta}\setminus N_{\alpha}$. First
note that $m\leqslant\frac{1}{\alpha-\beta}<m+1$ and therefore
$$m\beta<m\alpha\leqslant m\beta+1<m\beta+\beta=(m+1)\beta<\lfloor
(m+1)\beta\rfloor+1<(m+1)\alpha.$$
We distinguish two cases.

{\bf Case (1)}. $\lfloor m\beta\rfloor=\lfloor m\alpha\rfloor$. In
this case we have
\[\lfloor m\alpha\rfloor=\lfloor m\beta\rfloor<\lfloor
(m+1)\beta\rfloor<\lfloor (m+1)\beta\rfloor+1\leqslant\lfloor
(m+1)\alpha\rfloor.\]

\noindent This proves the claim.

{\bf Case(2)}. $\lfloor m\beta\rfloor\neq\lfloor m\alpha\rfloor$.
In this case we have $\lfloor m\alpha\rfloor=\lfloor
m\beta\rfloor+1$. Therefore \[\lfloor m\alpha\rfloor=\lfloor
m\beta\rfloor+1<(m\beta+1)+1\leqslant (m+1)\beta<\lfloor
(m+1)\beta\rfloor+1\leqslant (m+1)\beta+1<(m+1)\alpha.\] Hence
$\lfloor(m+1)\beta\rfloor<\lfloor
(m+1)\beta\rfloor+1\leqslant\lfloor (m+1)\alpha\rfloor$.

We show that $\lfloor m\alpha\rfloor<\lfloor (m+1)\beta\rfloor$.
Clearly $\lfloor m\alpha\rfloor\leqslant\lfloor
(m+1)\beta\rfloor$. Suppose $\lfloor m\alpha\rfloor=\lfloor
(m+1)\beta\rfloor$. Then $\lfloor m\beta\rfloor+1=\lfloor
(m+1)\beta\rfloor$. On the other hand, we have
\[\lfloor
m\beta\rfloor+2\leqslant m\beta+2\leqslant
m\beta+\beta=(m+1)\beta.\]

\noindent Thus $\lfloor (m+1)\beta\rfloor\geqslant\lfloor
m\beta\rfloor+2$, a contradiction. Therefore $\lfloor
m\alpha\rfloor=\lfloor m\beta\rfloor+1<\lfloor
(m+1)\beta\rfloor<\lfloor (m+1)\alpha\rfloor$ which again shows
the claim.
\end{proof}

\noindent Now we are going to study the above property in a general
$\langle F,I\rangle$.

\begin{definition} {\em We say that $\langle F,I\rangle$ is {\it
separable}, or $\langle F,I\rangle$ satisfies $\mathbb{S}$ property
for short, if for every distinct $\alpha,\beta\geqslant 1$,
$N_{\alpha}\neq N_{\beta}$ if and only if $\alpha\neq\beta$.}
\end{definition}

In Archimedean case, we prove $\mathbb{S}$ property by {\it
induction}. In fact we show that if
$\mathbb{N}_\alpha=\mathbb{N}_\beta$, then $\lfloor
n\alpha\rfloor=\lfloor n_\beta\rfloor$ for all $n\in\mathbb{N}$.
In non-Archimedean case, induction is too weak to can prove
$\mathbb{S}$. M. Moniri has claimed that ``{\it any arbitrary
structure $\langle F,I\rangle$ is separable}", (private
communication). We will show that a {\it weak} version of
$\mathbb{S}$ property can be deduced in every $\langle
F,I\rangle$. for this purpose, we need some auxiliary results. In
theorem \ref{rational} (ii), we showed that $N_{\alpha}\cap
N_{\beta}\neq\emptyset$ and $N_{\alpha}\cup N_{\beta}\neq
I^{\geqslant 0}$ for rational elements. When $\alpha$ and $\beta$
are irrational, the following result of Beatty \cite{Be}, also
reported in Skolem \cite{Sk}, provides a different view.

\begin{theorem} \label{Basic} Let $\alpha,\beta$ be positive
irrationals such that $\alpha^{-1}+\beta^{-1}=1$. Then
$N_{\alpha}\cap N_{\beta}=\{0\}$ and $ N_{\alpha} \cup
N_{\beta}=I^{\geqslant 0}$.
\end{theorem}

\begin{proof} To show $N_{\alpha}\cap N_{\beta}=\{0\}$,
suppose there exists $0\neq k\in N_{\alpha}\cap N_{\beta}$. Then
there would be $m,n\in I^{>0}$ such that $k\leqslant m\alpha<m+1,\
k\leqslant n\beta<k+1$. Since $\alpha,\beta$ are irrationals, so
the previous inequalities are proper. So
\[\frac{k}{m}<\alpha<\frac{k+1}{m},\ \frac{k}{n}<\beta<\frac{k+1}{n}.
\hspace{1 in}(1)\]

Thus $\frac{k}{n}<\frac{\alpha}{\alpha-1}<\frac{k+1}{n}$ where $
\beta=\frac{\alpha}{\alpha-1}$. Hence $
k(\alpha-1)<n\alpha<(k+1)(\alpha-1)$. But $\alpha,\beta>1$, so
$m,n<k$ and we have $(k-n)\alpha<k$ and $ k+1<(k+1-n)\alpha$. So
\[\frac{k+1}{k+1-n}<\alpha<\frac{k}{k-n}.\hspace{1 in} (2)\]

By (1) and (2), we have $\frac{k}{m}<\frac{k}{k-n}$ and $
\frac{k+1}{m}>\frac{k+1}{k+1-n}$, so $k-n<m<k+1-n$ implying
$k<m+n<k+1$. This would be a contradiction, since $m,n,k\in
I^{>0}$.

Next, we show $N_{\alpha}\cup N_{\beta}=I^{\geqslant 0}$. Suppose
there is some $h\in I^{>0}\setminus(N_{\alpha}\cap N_{\beta})$.
Then there exist $m,n\in I^{\geqslant 0}$ such that $\lfloor
m\alpha\rfloor<h<\lfloor (m+1)\alpha\rfloor$ and $\lfloor
n\beta\rfloor <h<\lfloor (n+1)\beta\rfloor$ implying
\[\begin{array}{lll}
m\alpha &<h &<(m+1)\alpha-1, \\
n\beta &<h &<(n+1)\beta-1. \\
\end{array} \]
From these two we get $(h+1)\alpha^{-1}-1<m<h\alpha^{-1}$ and $
(h+1)\beta^{-1}-1<n<h\beta^{-1}$. Therefore $
(h+1)(\alpha^{-1}+\beta^{-1})-2<m+n<h(\alpha^{-1}+\beta^{-1})$ and
so $h+1-2<m+n<h$ showing $h-1<m+n<h$. Since $m,n\in I^{\geqslant
0}$, the last inequality is impossible.
\end{proof}

\noindent We presented one direction of Fact A for real field in
Theorem \ref{Basic}. For real case, the proof of the converse of
Theorem \ref{Basic} is based on {\bf PHP} and some properties of
an auxiliary function. We define $\mu (\alpha,h)=|\{ n\in
\mathbb{N}\ | \ \lfloor n\alpha\rfloor\leqslant h\}|$, the number
of elements of $\mathbb{N}_\alpha$ not exceeding $h$. Note that in
the real case, if
$\mathbb{N}_\alpha\cup\mathbb{N}_\beta=\mathbb{N}$ and
$\mathbb{N}_\alpha\cap\mathbb{N}_\beta=\emptyset$, then
$$\mu(\alpha,h)+\mu(\beta,h)=h.\hspace{.5 in} (\ast)$$

\noindent The $(\ast)$ equality provides the proof of the converse
of Theorem \ref{Basic} in the field of real numbers. We have
$\lfloor\mu(\alpha,h)\alpha\rfloor\leqslant
h<\lfloor(\mu(\alpha,h)+1)\alpha\rfloor-1$. Using this
implication, we can define $\mu(\alpha,h)$ in non-Archimedean
case. But this definition doesn't provide $(\ast)$ and
consequently doesn't prove Fact A. On the other hand, if $\langle
F,I\rangle$ has the $\mathbb{S}$ property, we can deduce Fact A as
follows. Therefore study of the $\mathbb{S}$ property is a useful
tool for extend the theory of Diophantine approximations to
arbitrary ordered fields which equipped with integer parts.

\begin{lemma} Let $\langle F,I\rangle$ be separable. If
$N_{\alpha}\cap N_{\beta}={0}$ and $N_{\alpha} \cup
N_{\beta}=I^{\geqslant 0}$, then $\alpha,\beta$ are irrationals
and $\alpha^{-1}+\beta^{-1}=1$.
\end{lemma}

\begin{proof} By theorem \ref{rational}, one of $\alpha$, $\beta$ is
irrational. Suppose $\alpha$ is irrational. Set
$\alpha^{-1}+\eta^{-1}=1$. Then $N_{\eta}=N_{\beta}$ and
$\eta,\beta>1$. Since $\langle F,I\rangle$ is separable,
$\eta=\beta$.
\end{proof}

\noindent The following lemma is the rational version of Theorem
\ref{Basic}. Its proof is exactly similar to theorem \ref{Basic}.
The reader is refered to Theorem 3.15 in \cite{N}.

\begin{lemma} Let $\rho,\sigma$ be two positive rationals such
that $\rho^{-1}+\sigma^{-1}=1$. Then
\[ \begin{array}{ll}
N_{\rho}\cap N_{\sigma} &=\displaystyle{\bigcup_{m\in\rho I\cap
I^{\geqslant
0}}}mI^{\geqslant 0} \\
I^{\geqslant 0}\setminus(N_{\rho} \cup N_{\sigma})
&=\displaystyle{\bigcup_{m\in\rho I\cap I^{\geqslant
0}}}(mI^{\geqslant
0}+(m-1))\\
\end{array}\]
\end{lemma}

Note that $\rho I\cap I^{\geqslant 0}=\sigma I\cap I^{\geqslant
0}$. If $\rho=\frac{k}{m}$ is such that $(m,k)=1$, then
$\sigma=\frac{k}{k-m}$, and we have $N_{\rho}\cap
N_{\sigma}=kI^{\geqslant 0}$, and $I^{\geqslant
0}\setminus(N_{\rho} \cup N_{\sigma})=kI^{\geqslant 0}+(k-1)$. If
$\rho$ has no irreducible representation, then $\rho I^{\geqslant
0}\cap I^{\geqslant 0}$ is union of arithmetical progressions of
the form $m\rho I^{\geqslant 0}$ for some $m\in I^{> 0}$ such that
$m\rho\in I^{\geqslant 0}$. These arithmetical progressions have
nonempty intersections, but it is impossible to find $k\in
I^{\geqslant 0}$ such that $N_{\rho}\cap N_{\sigma}=kI^{\geqslant
0}$. The existence of an element like $``k"$ is equivalent to the
existence of an irreducible representation for $\rho$ as a
rational.

Now we can prove versions of the $\mathbb{S}$ property. The
following theorem proves this property for some large classes.

\begin{theorem} Let $\alpha,\beta>1$. Then $N_{\alpha}\neq
N_{\beta}$ when
\begin{enumerate}
\item $\alpha,\beta$ are distinct irrationals;
\item $\alpha,\beta$ are distinct rationals.
\end{enumerate}
\end{theorem}

\begin{proof}
First, suppose $\alpha$ and $\beta$ are distinct irrationals. It
suffices to prove the lemma for $1<\beta<\alpha<2$. There exist
$\eta,\gamma$ such that
$\eta^{-1}=1-\alpha^{-1},\gamma^{-1}=1-\beta^{-1}$ and
$2<\eta<\gamma$. So $N_{\gamma}\neq N_{\eta}$ and there exists $x\in
N_\eta\setminus N_\gamma$. Thus $x\in N_\beta\setminus N_{\alpha}$.

Now suppose $\rho,\sigma>1$ are distinct rationals. Using Theorem
\ref{T2}, we have only to consider the case $\rho,\sigma<2$. Let
$\eta$ and $\gamma$ be rational elements such that
$\rho^{-1}+\eta^{-1}=1$ and $\sigma^{-1}+\gamma^{-1}=1$. Then
$2<\gamma<\eta$ and by the above lemma, we have
\[ N_{\rho}=\left[ I^{\geqslant0}\setminus\left(N_{\eta}\cup
\displaystyle{\bigcup_{m\in\rho I\cap I^{\geqslant
0}}}(mI^{\geqslant 0}+(m-1))\right)\right]\cup(\rho I^{\geqslant
0}\cap I^{\geqslant 0}), \]

and
\[ N_{\sigma}=\left[ I^{\geqslant0}\setminus\left(N_{\gamma}\cup
\displaystyle{\bigcup_{m\in\sigma I\cap I^{\geqslant
0}}}(mI^{\geqslant 0}+(m-1))\right)\right]\cup(\sigma I^{\geqslant
0}\cap I^{\geqslant 0}). \]

\noindent Consequently by theorem \ref{T2}, there exists $x\in
N_{\gamma}\setminus N_{\eta}$. Let $m\in (\gamma I^{\geqslant
0}\cap\eta I^{\geqslant 0}\cap I^{\geqslant 0})$ such that
$x<m-1$. So $mI^{\geqslant 0}+x\subset N_{\gamma}\setminus
N_{\eta}$. Then $mI^{\geqslant 0}+x\subset N_{\rho}\setminus
N_{\sigma}$.
\end{proof}

\noindent By the same method as used in Theorem \ref{Basic}, we have

\begin{theorem} \label{irrational} Let $\alpha,\beta> 1$ be
distinct irrationals and $a,b,c,d\in I^{\geqslant 0}$. Then the
following properties hold
\begin{enumerate}
\item If $N_{\alpha}\cap N_{\beta}=\{0\}$
and $N_{\alpha} \cup N_{\beta}=I^{\geqslant 0}$, then
$\alpha^{-1}+\beta^{-1}=1$.

\item If $a\alpha^{-1}+b\beta^{-1}=1$, then $N_{\alpha}\cap N_{\beta}=\{0\}$.

\item If $a(1-\alpha^{-1})+b(1-\beta^{-1})=1$, then
$N_{\alpha} \cup N_{\beta}=I^{\geqslant 0}$.

\item If  $a\alpha^{-1}+b(1-\beta^{-1})=1$, then $
N_{\alpha}\subseteq N_{\beta}$.

\item If
$a\alpha^{-1}+b\beta^{-1}=1,c(1-\alpha^{-1})+d(1-\beta^{-1})=1$,
then $a=b=c=d=1$ (and so $\alpha^{-1}+\beta^{-1}=1$).
\end{enumerate}
\end{theorem}

\begin{proof} {\bf 1.} This is obvious, since $\alpha$ and
$\beta$ are irrational.

{\bf 2.} Suppose there exists $0\neq k\in (N_{\alpha}\cap
N_{\beta})$. Then there exist $m,n\in I^{>0}$ such that
$k<m\alpha<k+1$ and $k<n\beta<k+1$. Consequently
$$\frac{k}{m}<\alpha<\frac{k+1}{m} and
\frac{k}{n}<\beta<\frac{k+1}{n}\ {\bf (1)}.$$ Therefore
$\frac{k}{m}<\frac{b\alpha}{\alpha-a}<\frac{k+1}{n}$, where $
\beta=\frac{b\alpha}{\alpha-a} $. Hence
$k(\alpha-a)<nb\alpha<(k+1)(\alpha-a)$ and so $(k-bn)\alpha<ak$
and $(k+1)a<(k+1-bn)\alpha$. We claim that $k-bn>0$. To see this,
suppose $\beta>\alpha$. Then $a\alpha^{-1}+b\beta^{-1}=1$.
Therefore $an\beta\alpha^{-1}+bn=n\beta$. So
$an\beta\alpha^{-1}+bn<k+1$ and therefore $k+1-bn>1$ implying
$k-bn>0$. Therefore we have
$\frac{a(k+1)}{k+1-bn}<\alpha<\frac{ak}{k-bn}$. Using {\bf (1)},
we have $\frac{k}{m}<\frac{ak}{k-bn}$ and
$\frac{a(k+1)}{k+1-bn}<\frac{k+1}{m}$ and so $k<am-bn<k+1$. But
$am-bn\in I$, a contradiction.

{\bf 3.} Let $\gamma^{-1}=1-\alpha^{-1},\eta^{-1}=1-\beta^{-1}$.
Then $\eta$ and $\gamma$ are irrationals and
$a\gamma^{-1}+b\eta^{-1}=1$. Therefore $N_{\gamma}\cap\
N_{\eta}={0}$. Since we have $N_{\gamma}\cap\ N_{\alpha}={0}$ and
$N_{\gamma}\cup\ N_{\alpha}=I^{\geqslant 0}$, so
$N_{\eta}\subseteq N_{\alpha}$ and finally $N_{\beta}\cup
N_{\alpha}\supseteq N_{\beta}\cup N_{\eta}=I^{\geqslant 0}$.

{\bf 4.} Let $\eta^{-1}=1-\beta^{-1}$. Then
$a\alpha^{-1}+b\eta^{-1}=1$. Therefore $N_{\eta}\cap
N_{\alpha}={0}$ and so $N_{\alpha}\subseteq N_{\beta}$.

{\bf 5.} We have $N_{\alpha}\cap N_{\beta}={0}$ and
$N_{\alpha}\cup N_{\beta}=I^{\geqslant 0}$. Since $\alpha,\beta$
are irrationals, $\alpha^{-1}+\beta^{-1}=1$. Therefore
$a\alpha^{-1}+b\beta^{-1}=\alpha^{-1}+\beta^{-1}$. Hence we have
$\alpha^{-1}(a-1)=\beta^{-1}(1-b)\geqslant 0$. From $b\geqslant
1$, we find that $b=a=1$. Also, $c=d=1$.
\end{proof}

We now must determine the relation between $N_{\alpha}$ and
$N_{\rho}$, when $\alpha$ is irrational and $\rho$ is rational.
Note that $Q$, the fraction field of $I$ is a dense subfield of
$F$ and if $F$ has an irrational element, $Q$ and $F\setminus Q$
are proper dense subsets of $F$. So, if $\rho$ is a rational
element, then for each positive $\epsilon\in F$, there exist some
irrationals $\alpha$ such that $|\alpha-\rho|<\epsilon$. Using
this property, it is easy to define convergent sequences in scale
of the ordered field $F$. Therefore if $cf(I)=\eta$, we have some
$\eta$-sequences of irrationals which converge to $\rho$, (Note
that $cf(F)=cf(I)$). The following considers this situation.

\begin{theorem} Let $\rho\geqslant 1$ be a rational. Suppose
$cf(I)=\eta$, and $\{\alpha_{\gamma}\}_{\gamma <\eta}$ is a
descending sequence of irrationals such that
$\lim_{\gamma\rightarrow\eta}\alpha_{\gamma}=\rho$. Then for every
$m\in I^{> 0}$, there exists $\beta<\eta$ such that for all
$\beta<\gamma<\eta$, we have
$N_{\rho}|_{<m}=N_{\alpha_{\gamma}}|_{<m}$.
\end{theorem}

\begin{proof} Suppose $\rho=\frac{p}{q}$ is such that $p,q\in
I^{>0}$. We can assume that $m=qt$. Otherwise, consider a multiple
of $q$ greater than $m$, such as $(\lfloor\frac{m}{q}\rfloor+1)q$.
Then $\rho=\frac{pt}{qt}=\frac{pt}{m}$ and therefore, for all
$l\in I^{\geqslant 0}$ if $l\leqslant m$, there exists $u_l\in
I^{\geqslant 0}$ such that $\rho\in
[k+\frac{u_l}{l},k+\frac{u_l+1}{l})$, where $k= \lfloor \rho
\rfloor$. We show that there exists an interval which contains
$\rho$ in the intersection of $\bigcap_{l\leqslant m}
[k+\frac{u_l}{l},k+\frac{u_l+1}{l})$. For each interval
$[k+\frac{u_l}{l},k+\frac{u_l+1}{l})$, we have two cases:
($\rho=k+\frac{u_m}{m}$, $l<m$)

{\bf Case (1).} $k+\frac{u_l}{l}\in
[k+\frac{u_m}{m},k+\frac{u_m+1}{m})$. In this case, we have
$\rho=k+\frac{u_l}{l}$.

{\bf Case (2).} $k+\frac{u_l+1}{l}\in
[k+\frac{u_m}{m},k+\frac{u_m+1}{m}]$. In fact, the two cases are
disjoint, since the length of $[k+\frac{u_l}{l},k+\frac{u_l+1}{l})$
is $\frac{1}{l}>\frac{1}{m}$. In case(2), we have
\[ k+\frac{u_l+1}{l}-\rho=k+\frac{u_l+1}{l}-(k+\frac{u_m}{m})
=\frac{u_l+1}{l}-\frac{u_m}{m}\geqslant\frac{1}{ml}>\frac{1}{m^2}.\]
Therefore there exists an interval $I$ with length at least
$\frac{1}{m^2}$ such that for all $l$ less than $m$,
$I\subseteq[k+\frac{u_l}{l},k+\frac{u_l+1}{l})$ contains $\rho$.

On the other hand, $\alpha_{\gamma}\searrow\rho$ and the
$\alpha_{\gamma}$'s are irrationals. So there exists $\beta<\eta$
such that for all $\gamma>\beta$, we have
$0<\alpha_{\gamma}-\rho<\frac{1}{m^2}$ and therefore, for all
$l\leqslant m$, $\alpha_{\gamma}\in
[k+\frac{u_l}{l},k+\frac{u_l+1}{l})$. So for all $\gamma>\beta$, $
\alpha_{\gamma}\in I$. But for all $x\in I$ (such as $\rho$ and
$\alpha_{\gamma}$ for $\gamma>\beta$) and $l\leqslant m$, we have
$k+\frac{u_l}{l}\leqslant x<k+\frac{u_l+1}{l}$ and so
$kl+u_l\leqslant lx<kl+u_l+1$ implying $\lfloor lx\rfloor=kl+u_l$.
Since $km+u_m\geqslant m$, we get $\beta<\gamma<\eta$ and so
$N_{\rho}|_{<m}=N_{\alpha_{\gamma}}|_{<m}$.
\end{proof}

Below, using the above theorem, we generalize one direction of
Fact F' without any condition and show that if $I$ is a B\'ezout
domain, Fact F' generalizes.

\begin{theorem} Suppose $\rho,\sigma>1$ are rationals. Then we have
\begin{enumerate}
\item if there exist $a,b\in I^{\geqslant 0}$ such that $\
a\rho^{-1}+b(1-\sigma^{-1})=1$. Then $N_{\rho} \subseteq
N_{\sigma}$.
\item if $I$ is B\'ezout and $N_{\rho}\subseteq
N_{\sigma}$, then there exist $a,b\in I^{\geqslant 0}$ such that
$a\rho^{-1}+b(1-\sigma^{-1})=1$.
\end{enumerate}
\end{theorem}

\begin{proof}
{\bf 1}. Let $\{\alpha_{\eta}\}$ be a decreasing sequence of
irrationals which tends to $\sigma$. For every sufficiently large
ordinal $\eta$, $N_{\alpha_{\eta}}$ and $N_{\sigma}$ coincide on
some initial segment with arbitrary large length. Choose
$\{\beta_{\eta}\}$ such that
$a(1-\alpha_{\eta})^{-1}+b\beta_{\eta}^{-1}=1$. Then the sequence
$\{\beta_{\eta}\}$ is a decreasing sequence converging to $\rho$.
Now suppose $t\in N_{\rho}$. Then there is an ordinal $\gamma$ such
that for all $\eta>\gamma$, $N_{\beta_{\eta}}$ and $N_{\rho}$
coincide on $\leqslant t$. Note that $t\in N_{\beta_{\eta}}$. Since
$N_{\beta_{\eta}}\subseteq N_{\alpha_{\eta}}$, we must have $t\in
N_{\alpha_{\eta}}$. In particular, $t\in N_{\alpha}$.

{\bf 2}. We show that if $\sigma\in I^{>1}$,  then $\rho\in
I^{>0}$ and consequently $\rho$ is a multiple of $\sigma$. Take
$\rho=\frac{m}{n}$ with $(m,n)=1$, $m>n>1$ . Since $I$ is
B\'ezout, there exist $s,t\in I$ such that $sm+tn=1$. If $s<0$,
take $k\in I^{\geqslant 0}$ such that $k>[\frac{s}{n}]$. Then
$s+kn\in I^{\geqslant 0}$ and $m(s+kn)+n(t-km)=1$. Suppose $s>0$.
Then $ms=1+nk'$, $[s\rho]=[s\times \frac{m}{n}]=k'$,
$[ns\rho]=ms$, and $(ms,k')=1$. But $ms,k'$ are both multiples of
$\sigma\in I^{>0}$, a contradiction.

Now suppose $\sigma\in Q\setminus I$, $\sigma=\frac{m}{s}$ and let
$\rho=\frac{m}{n}$ be such that $(m,n,s)=1$. If $(m-s,n)=d$, then
$\frac{m}{d}$ will be an integer multiple of both $\frac{m}{n}$
and $\frac{m}{m-s}$ (i.e., there will exist $u,v\in I$ such that
$u\cdot \frac{m}{n}=\frac{m}{d}$ and $v\cdot
\frac{m}{m-s}=\frac{m}{d}$). Then $N_{\frac{m}{d}}\subseteq
N_{\frac{m}{m-s}}$ and $N_{\frac{m}{d}}\subseteq
N_{\frac{m}{n}}\subseteq N_{\frac{m}{s}}$. So,
$N_{\frac{m}{d}}\subseteq N_{\frac{m}{s}}\cap N_{\frac{m}{m-s}}$.
But if $\eta=\frac{m}{m-s}$, then $\eta^{-1}+\sigma^{-1}=1$ and if
$k=\frac{m}{(m,s)}$, then $N_{\frac{m}{s}}\cap N_{\frac{m}{m-s}}=k
I^{\geqslant 0}$. Therefore $N_{\frac{m}{d}}\subseteq N_{k}$. By
the first paragraph, $d\mid m$ and $d\mid (m,n,s)=1$. So,
$(m-s,n)=(m,n,s)=1$.

For all $j\in I^{\geqslant 0}$, there exists $x_j\in I^{\geqslant
0}$ such that $[j\rho]=[x_j\sigma]$. Thus $j\rho-y,xj\sigma-y$
have the same signs for all $y\in I$ (the sign of zero is taken
here to be plus). Substituting $\rho$ and $\sigma$, we conclude
that $jm-yn$ and $x_jm-ys$ have the same signs for all $j$ and
$y$. We use absolute values of these numbers as $a$ and $b$, with
$j$ and $y$ chosen appropriately. First note that
$a(1-\sigma^{-1}+b\rho^{-1})=\mid(jm-yn)(1-\frac{m}{s})+(x_jm-ys)\frac{n}{m}\mid=\mid
j(m-s)+(x_j-y)n\mid$. For any fixed positive integer $j$, note
that $x_j-y$ can assume all $I$-values. This argument proves the
theorem unless one or the other of these values $a$ or $b$ is $0$.
But as the standard case (i.e., in the real field), we can choose
appropriate $a$'s and $b$'s.
\end{proof}

\begin{corollary}Suppose $\rho,\sigma>1$ are rationals and
$I$ is B\'ezout. Then $N_{\rho}\subseteq N_{\sigma}$ if and only
if there exist $a,b\in I^{\geqslant 0}$ such that
$a\rho^{-1}+b(1-\sigma^{-1})=1$.
\end{corollary}

\noindent Now we show that if $N_\sigma=N_\rho$, then
$\sigma,\rho$ are close to each other.
\begin{lemma}\label{Linf} Let $\sigma,\rho\in F^{\geqslant 1}$ with $\sigma
<\rho$ be such that $N_{\sigma}=N_{\rho}$. Then $\rho-\sigma$ is
an infinitesimal.
\end{lemma}

\begin{proof} We already know that $1\leqslant\sigma<\rho<2$ and
so $\lfloor\sigma\rfloor=\lfloor\rho \rfloor=1$. Assume for the
sake of a contradiction that $\rho-\sigma$ is not an infinitesimal
and so its inverse is limited. Suppose $\rho-\sigma=r+\epsilon$
such that $r$ is a real number and $\epsilon$ is an infinitesimal.
Then $0\leqslant r\leqslant 1$, (note that if $r=1$ then
$\epsilon<0$). Thus $\frac{1}{\rho-\sigma}=\frac{1}{r+\epsilon}$
and
$\frac{1}{r}-\frac{1}{r+\epsilon}=\frac{\epsilon}{r(r+\epsilon)}$
is an infinitesimal because $r$ and $r+\epsilon$ are both finite.
Therefore $\frac{1}{\rho-\sigma}$ has a standard integer part
$m=\lfloor\frac{1}{\rho-\sigma}\rfloor$ which equals to
$\lfloor\frac{1}{r}\rfloor$ or $\lfloor\frac{1}{r}\rfloor-1$.
Therefore $m$ is {\bf finite}. We have $(m+1)\sigma+1<(m+1)\rho$,
since $\frac{1}{\rho-\sigma}<m+1$. Thus
$\lfloor(m+1)\sigma\rfloor<\lfloor(m+1)\rho\rfloor$. So there
exists $1\leqslant k\leqslant m$ such that $\lfloor k
\rho\rfloor=\lfloor k\sigma \rfloor$ and
$\lfloor(k+1)\rho\rfloor>\lfloor (k+1)\sigma \rfloor$. Thus
$\lfloor (k+1) \sigma \rfloor\in N_{\sigma}\setminus N_{\rho}.$
\end{proof}

We prove that if $2<\alpha<\beta$, then
$\lfloor(m+1)\alpha\rfloor\in N_{\alpha}\setminus N_{\beta}$ for
$m=\lfloor\frac{1}{\beta-\alpha}\rfloor$. Now let
$1<\alpha<\beta<2$ and both are irrationals. So if
$\alpha^{-1}+\eta^{-1}=1$ and $\beta^{-1}+\gamma^{-1}=1$, then
$2<\gamma<\eta$ and we have $\lfloor(m+1)\alpha\rfloor\in
N_{\gamma}\setminus N_{\eta}$ for
$m=\lfloor\frac{1}{\eta-\gamma}\rfloor$ and consequently it is in
$N_{\alpha}\setminus N_{\beta}$. Note that $m=\lfloor
\frac{(\alpha-1)(\beta-1)}{\beta-\alpha}\rfloor$. We use from this
element to prove the $\mathbb{S}$ property for a suitable case.

Suppose $\rho<\beta$ and $\rho$ is rational and $\beta$ is
irrational. Then for all irrationals sufficiently close to $\rho$
(and greater than it), such as $\alpha$, if
$m=\lfloor\frac{(\alpha-1)(\beta-1)}{\beta-\alpha}\rfloor$, then
we have $\lfloor(m+1)\frac{\beta}{\beta-1}\rfloor\in\
N_{\frac{\beta}{\beta-1}}\setminus N_{\frac{\alpha}{\alpha-1}}$.
So, $\lfloor(m+1)\frac{\beta}{\beta-1}\rfloor\in
N_{\alpha}\setminus N_{\beta}$ and thus
\[\frac{(\rho-1)(\beta-1)}{\beta-\rho}-\frac{(\alpha-1)(\beta-1)}
{\beta-\alpha}=\frac{(\rho-1)(\alpha-\rho)}{(\beta-\rho)(\beta-\alpha)}(1-\beta)<0.\]
So
$\frac{(\rho-1)(\beta-1)}{\beta-\rho}<\frac{(\alpha-1)(\beta-1)}
{\beta-\alpha}$ and we have
$\lfloor(m+1)\frac{\rho}{\rho-1}\rfloor=\lfloor(m+1)\frac{\alpha}{\alpha-1}\rfloor$
for all irrationals $\alpha$ which are sufficiently closed to
$\rho$ (and greater than it). So
$\lfloor(m+1)\frac{\alpha}{\alpha-1}\rfloor\in N_{\rho}\setminus
N_{\beta}$.

If $\alpha<\rho<\beta$ and $\rho$ is rational. We have
$\frac{(\alpha-1)(\rho-1)}{\rho-\alpha}-\frac{(\alpha-1)(\beta-1)}
{\beta-\alpha}=\frac{(\alpha-1)(\rho-\beta)}{(\rho-\alpha)(\beta-\alpha)}<0$.
So
$\frac{(\alpha-1)(\rho-1)}{\rho-\alpha}<\frac{(\alpha-1)(\beta-1)}
{\beta-\alpha}$. Fix $1<\alpha<\rho<2$, such that $\alpha$ is
irrational and $\rho$ is rational. If
$m=\lfloor\frac{(\alpha-1)(\rho-1)}{\rho-\alpha}\rfloor$, then for
all irrationals sufficiently close to $\rho$ (and larger than it),
called it $\beta$, we imply that
$\lfloor(m+1)\frac{\beta}{\beta-1}\rfloor\in N_{\alpha}\setminus
N_{\beta}$. We have
$\frac{\rho}{\rho-1}-\frac{\beta}{\beta-1}=\frac{\beta-\rho}{(\beta-1)(\rho-1)}>0$.
So if $\beta$ is sufficiently close to $\rho$, then $\lfloor (m+1)
\frac{\rho}{\rho-1}\rfloor=\lfloor (m+1)\frac{\beta}{\beta-1}
\rfloor$, unless $(m+1)\frac{\rho}{\rho-1}\in I^{>0}$. So,
$\lfloor (m+1)\frac{\rho}{\rho-1}\rfloor\in N_{\alpha}\setminus
N_{\rho}$, unless $(m+1)\frac{\rho}{\rho-1}\in I^{>0}$. If
$(m+1)\frac{\rho}{\rho-1}\in I^{>0}$, then
$(m+1)\frac{\rho}{\rho-1}\in N_{\alpha}\cap N_{\rho}\cap
N_{\frac{\rho}{\rho-1}} \cap I^{>0}$.

Totally, if $\rho$ is rational and $\beta$ is irrational, we prove
the $\mathbb{S}$ Property for all $1<\rho<\beta<2$ and for all
$1<\beta<\rho<2$ s.t. $ (m+1) \frac{\rho}{\rho-1} \not\in I$ for
$m=\lfloor\frac{(\beta-1)(\rho-1)}{\rho-\beta}\rfloor$. An example
for the last case, in $\langle\mathbb{R},\mathbb{Z}\rangle$, let
$\rho=\frac{3}{2}$ and $\beta=\sqrt{2}$. Then $m=2$,
$(m+1)\frac{\rho}{\rho-1}=9$, $7\times\beta=\lfloor
7\sqrt{2}\rfloor=9$ and $6\times\rho=9$.

Suppose $1<\alpha<\beta<2$ and $\alpha,\beta$ are two arbitrary
elements of $F$ such that $N_{\alpha}=N_{\beta}$. Then one of them
is rational and the other is irrational. So for all $\gamma\in F$
which are $\alpha<\gamma<\beta$, $N_{\gamma}\neq N_{\beta}$.
Because:
\begin{enumerate}
\item Let $\alpha$ be irrational and $\beta$ be
rational, if $\gamma$ is irrational, $N_\gamma\neq N_\alpha$ and
if $\gamma$ is rational, $N_\gamma\neq N_\beta$.

\item Let $\alpha$ be rational and $\beta$ be irrational,
if $\gamma$ is irrational, $N_\gamma\neq N_\beta$ and if $\gamma$
is rational, $N_\gamma\neq N_\alpha$.
\end{enumerate}

\noindent Now we define $\alpha\thicksim\beta$ if
$N_\alpha=N_\beta$. This relation is an {\it equivalence relation}
and we have:

\noindent (1) if $0<\alpha\leqslant 1$, then
$[\alpha]_\thicksim=(0,1]_F$;

\noindent (2) if $\alpha\geqslant 2$, then
$[\alpha]_\thicksim=\{\alpha\}$;

\noindent (3) If $\langle F,I\rangle$ is separable and $\alpha>1$
then $[\alpha]_\thicksim=\{\alpha\}$;

\noindent (4) for an arbitrary $\langle F,I\rangle$, if $\alpha>1$
then $[\alpha]_\thicksim=\{\alpha\}$ or
$[\alpha]_\thicksim=\{\alpha,\beta\}$ such that if $\alpha$ is
irrational then $\beta$ is rational and vise versa and
$\alpha-\beta$ is an infinitesimal element in $F$.
\section{Arithmetical Progressions }

In real case, for any irrational $\alpha>0$, the set
$\mathbb{N}_{\alpha}$ has a number of interesting number theoretic
properties. For example, for each $k<m\in\mathbb{N}$, the subset
$\{x|\ x\in\mathbb{N}_{\alpha},\ x\equiv k\ (mod\ m)\}$ is
unbounded and $\mathbb{N}_{\alpha}$ is uniformly distributed
modulo of every $m\in \mathbb{N}$, \cite{N}. In this section, we
will show that the first property is always equivalent to the
$\DMO$ property which is stronger than the $\mathbb{S}$ property.
\subsection{P Condition}

\begin{definition}A set $D\subset F$ is {\it dense modulo one} (or $\DMO$) with
respect to $I$ if the set $\{u-\lfloor u\rfloor_I|u\in D\}$ is dense
in $[0,1)_F$.
\end{definition}
In \cite{AM}, we presented some non-trivial $\DMO$ sets. Let recall
one of those example.

\begin{proposition} For every $\langle F,I\rangle$ and $p\in
\mathbb{N}$, the set $\{\sqrt[p]{u}\ |\ u\in I^{>0}\}$ is $\DMO$
with respect to every $\IP$ for $F$.
\end{proposition}

\begin{proof} Let $I_{1}$ be an $\IP$ for $F$. Suppose $k,t\in
I_1^{\geqslant 0}$ and $k<t$. We need to find $M\in I$ and $n\in
I_1$ such that $\frac{k}{t}<\sqrt[p]{M}-n<\frac{k+1}{t}$. We have
$n+\frac{k}{t}<\sqrt[p]{M}<n+\frac{k+1}{t}$, equivalently,
$(n+\frac{k}{t})^p<M<(n+\frac{k+1}{t})^p$. But
$(n+\frac{k+1}{t})^p-(n+\frac{k}{t})^p\in Frac(I_1)$ and we have
$$(n+\frac{k+1}{t})^p-(n+\frac{k}{t})^p=\frac{1}{t}((n+\frac{k+1}{t})^{p-1}
+(n+\frac{k+1}{t})^{p-2}(n+\frac{k}{t})\ldots
+(n+\frac{k}{t})^{p-1}).$$ Note that this is greater than
$\frac{p}{t}(n+\frac{k}{t})^{p-1}$. So if we choose $n\in I_1$ such
that the latter is greater than 1 (it suffices to choose
$\lfloor\sqrt[p-1]{\frac{t}{p}}\rfloor+1\leqslant n$), then there
will exist $M\in I^{>0}$ such that
$n+\frac{k}{t}<\sqrt[p]{M}<n+\frac{k+1}{t}$.
\end{proof}

\begin{definition}For an irrational $\alpha>0$, we say that $\DMO(\alpha)$ holds
whenever the set $D=\{n\alpha|n\in I^{\geqslant 0}\}$ is $\DMO$
(with respect to $I$).
\end{definition}
\begin{theorem} \label{Dmo} If $1\leqslant\alpha \in F$, the
following are equivalent:

(a) $\DMO(\alpha)$,

(b) $(\forall m,k\in I^{\geqslant 0})(N_{m\alpha}\cap(m\
I^{\geqslant 0}+k)\neq \{0\}).$
\end{theorem}

\begin{proof} $\bf {(a)\rightarrow (b)}$. By assumption, for all
$k,m\in I^{\geqslant 0}$, there exists $u\in I^{\geqslant 0}$ such
that $0\leqslant\frac{k}{m}\leqslant u\alpha-\lfloor
u\alpha\rfloor<\frac{k+1}{m}<1$. This shows $k\leqslant
mu\alpha-m\lfloor u\alpha\rfloor<k+1$ and so $m\lfloor
u\alpha\rfloor+k\leqslant mu\alpha<m\lfloor u\alpha\rfloor+k+1$
which in turn implies $\lfloor m u\alpha\rfloor\equiv k$, (mod m).

$\bf {(b)\rightarrow (a)}$. Suppose $0<l<r<1$. Since $Q=Frac(I)$
is a dense subfield of $F$, so there exist $p,q\in I^{\geqslant
0}$ such that $l<\frac{p}{q}<\frac{p+1}{q}<r$ (it suffices to
assume $\frac{1}{q}<r-l$). By $(b)$, there exists $n\in
I^{\geqslant 0}$ such that $\lfloor nq\alpha\rfloor\equiv p$ $(mod
q)$. Then $\exists t\in I^{\geqslant 0} qt+p\leqslant
nq\alpha<qt+p+1$ and so $t+\frac{p}{q}\leqslant
n\alpha<t+\frac{p+1}{q}$. This implies $l<\frac{p}{q}\leqslant
n\alpha-t<\frac{p+1}{q}<r$. Since $0<n\alpha-t<1$, $t=\lfloor
n\alpha\rfloor$ and therefore $l<n\alpha-\lfloor
n\alpha\rfloor<r$.
\end{proof}

In the field of real numbers, if $\alpha>1$ is irrational, then
$\mathbb{N}_{\alpha}$ intersects any arithmetical progression, but
does not contain any of them, (see \cite[Theorem 3.3]{N}).
Therefore, Theorem \ref{Dmo} gives another proof for
$\DMO(\alpha)$ in the standard situation. However we don't know
whether $\DMO(\alpha)$ holds in general or not. In the real case,
$\DMO(\alpha)$ for an irrational number $\alpha$ is usually
obtained via cofinal rational quadratic approximations. We deal
with this issue in Section 4. If $\DMO(\alpha)$ holds, then by
Theorem \ref{Dmo} and Theorem \ref{rational}, $\alpha>1$ will be
an irrational.

\begin{corollary} Suppose $\alpha>1$ is irrational. Then the two conditions in Theorem \ref{Dmo} are
equivalent to

(c) $(\forall m\in I^{> 0})(N_{m\alpha}\cap m\ I^{\geqslant 0}\neq
\{0\}).$
\end{corollary}

\begin{proof} The property $\DMO(\alpha)$ holds if and only if
for all $\epsilon
>0$, there exists some $n\in I^{\geqslant 0}$ such that
$n\alpha-\lfloor n\alpha\rfloor<\epsilon$. The reason goes as
follows. Pick $0<l<r<1$ and let $\epsilon=r-l$. There exists $n\in
I^{\geqslant 0}$ such that $n\alpha-\lfloor
n\alpha\rfloor<\epsilon$. Therefore $\frac{r-l}{n\alpha-\lfloor
n\alpha\rfloor}>1$. So there exists $k\in I^{>0}$ such that
$\frac{l}{n\alpha-\lfloor n\alpha\rfloor}
<k<\frac{r}{n\alpha-\lfloor n\alpha\rfloor}$. Thus we have
$0<l<kn\alpha-k\lfloor n\alpha\rfloor<r<1$ and therefore $k\lfloor
n\alpha\rfloor=\lfloor kn\alpha\rfloor$. Now, let $m=kn$, and so
$l<m\alpha-\lfloor m\alpha\rfloor<r$.

We have $b\rightarrow c$. Now let $\epsilon>0$. Set $m\in I$ such
that $m>\epsilon^{-1}$. Because of $(c)$, there exists $k,t\in
I^{>0}$, $mt<km\alpha<mt+1$. Therefore
$t<k\alpha<t+\frac{1}{m}<\epsilon$ and thus
$0<k\alpha-t<\epsilon$. Using the previous paragraph, proof is
complete.
\end{proof}

\noindent Now using the Theorem \ref{Dmo}, we present a new
property for structures $\langle F,I\rangle$:

\begin{definition} {\em Let $\alpha>1$ be an irrational. We
define four properties as follows:

$\mathbb{P}_{\alpha}^{1}$: The set $N_{\alpha}$ intersects each
arithmetical progressions,

$\mathbb{P}_{\alpha}^{2}$: The set $N_{\alpha}$ does not contain any
arithmetical progressions,

$\mathbb{P}$: For all irrationals $\alpha>1$,
$\mathbb{P}_{\alpha}^{1}$ holds,

$\mathbb{P}'$: For all irrationals $\alpha>1$,
$\mathbb{P}_{\alpha}^{2}$ holds.}
\end{definition}
\begin{proposition}\label{P1} If $\alpha,\beta>1$ are
irrationals such that $\beta^{-1}+\alpha^{-1}=1$. Then

(i) If $\mathbb{P}_{\alpha}^{2}$, then $(\forall m\in I^{>
0})(\mathbb{P}_{m\alpha}^{2})$.

(ii) We have $(\forall m\in I^{> 0})(\mathbb{P}_{m\alpha}^{1})$ if
and only if $\DMO(\alpha)$.

(iii) The properties $\mathbb{P}^{2}_{\alpha}$ and
$\mathbb{P}_{\beta}^{1}$ are equivalent.
\end{proposition}

\begin{proof} (i) It is sufficient to observe that for all $m\in
I^{\geqslant 0}$, we have $N_{m\alpha}\subseteq N_{\alpha}$.

(ii) This is just the content of Theorem \ref{Dmo}.

(iii) Suppose that $\mathbb{P}^2_{\alpha}$ holds. Then $N_\alpha$
does not contain any arithmetical progressions. So $N_\beta$ has a
nonempty intersection with every arithmetical progression by
Theorem \ref{Basic}.
\end{proof}

The structure $\langle F,I\rangle$ satisfies $\DMO$ if $F$ has
irrational elements and for all irrational element $\alpha$,
$\DMO(\alpha)$ hold. Using the similar method, we have the
following lemma.
\begin{lemma}\label{LP} Suppose that $\alpha>1$ is irrational. Then we
have the following:

(i) If $N_{\alpha}\cup N_{\gamma}=I^{\geqslant 0}$ and
$\mathbb{P}_{\alpha}^{2}$ holds, then so does
$\mathbb{P}_{\gamma}^{1}$.

(ii) If $N_{\alpha}\cap N_{\gamma}=\{0\}$ and
$\mathbb{P}_{\gamma}^{1}$ holds, then so does
$\mathbb{P}_{\alpha}^{2}$.

(iii) The $\mathbb{P}$ property holds if and only if $\mathbb{P}'$
does.

(iv) The $\mathbb{P}$ property holds if and only if $(\forall
\alpha>1)$ with $\alpha \in F\setminus Q$, we have $\DMO(\alpha)$
or more continently $\langle F,I\rangle\models\mathbb{P}$ if and
only if $\langle F,I\rangle\models\DMO$
\end{lemma}

In this section we study the structures $\langle F,I\rangle$ which
satisfies the $\mathbb{P}$ property. By part (iv) of Lemma
\ref{LP}, this section is about $\DMO$-$\langle F,I\rangle$, i.e.
the structures $\langle F,I\rangle\vDash\DMO$. At first, by using
the method similar to the above Lemma, we immediately get the
following theorem.

\begin{theorem}\label{P11} Suppose that $\langle F,I\rangle$ satisfies
the $\mathbb{P}$ property and $\alpha>1$ is irrational and
$\rho>1$ is rational. Then neither of the relations below could
hold:

$N_{\alpha}\subseteq N_{\rho}$, $N_{\rho}\subseteq N_{\alpha}$,
$N_{\alpha}\cap N_{\rho}=\{0\}$, $N_{\alpha}\cup
N_{\rho}=I^{\geqslant 0}$.
\end{theorem}

\begin{proof} Suppose $\rho=\frac{p}{q}$. Then
$N_\rho\cap(pI^{\geqslant 0}+(p-1))=\emptyset$. By the
$\mathbb{P}$ property, $N_\alpha$ has a nontrivial intersection
with this arithmetical progression. Hence $N_\alpha\nsubseteqq
N_\rho$.

The set $N_\rho$ is a union of arithmetical progressions and
$\mathbb{P}^2_\alpha$ holds. Therefore $N_\rho\nsubseteqq N_\alpha$.

By $\mathbb{P}^1_\alpha$, we have $N_\rho\cap N_\alpha\neq \{0\}$.

We have $(pI^{\geqslant 0}+(p-1))\setminus N_\alpha\neq\emptyset$
and $(pI^{\geqslant 0}+(p-1))\cap N_\rho=\emptyset$. So,
$N_\rho\cup N_\alpha\neq I^{\geqslant0}$.
\end{proof}

Using the above theorem, we introduce the relation between
separability and the $\mathbb{P}$ property.

\begin{corollary} If $\langle F,I\rangle$ has the $\mathbb{P}$ property,
then it is separable.
\end{corollary}

\begin{proof} We have shown already that $\langle F,I\rangle$ has
the $\mathbb{S}$ property if and only if for every
$1\leqslant\beta<\alpha<2$, with one of $\alpha,\beta$ being
rational and the other irrational, $N_{\alpha}\neq N_{\beta}$.
Theorem \ref{P11} completed the proof.
\end{proof}
\subsection{P Condition \& Skolem-Bang's Theorems}

It can be shown that the $\mathbb{P}$ property is a first order
sentence in $\langle F,I\rangle$. So by Upward L\"owenheim-Skolem
theorem over $\langle \mathbb{R},\mathbb{Z}\rangle$ (or over the
countable structure
$\langle\widetilde{\mathbb{Q}},\mathbb{Z}\rangle$), there exist
sufficiently large models of $\langle
F,I\rangle\models\mathbb{P}$. Professor Moniri conjectures the
following (private communication):

\begin{center} `` $\langle F,I\rangle\models\mathbb{P}$, for all ordered
field $F$ with IP $I$."
\end{center}

\noindent But now we want to discuss about the $\mathbb{P}$
property and Skolem-Bang's Theorems. In this subsection, suppose
$\langle F,I\rangle\models \mathbb{P}$.

\begin{theorem} If $I\models B\acute{e}z$ and $\alpha,\beta$ are
positive irrationals such that $1,\alpha,\beta$ are linearly
dependent over the $Frac(I)$, say

\begin{center} $a\alpha+b\beta=c$, $(a,b,c)=1$ and $c>1$,
\end{center}

\noindent then the points $(m\alpha-\lfloor
m\alpha\rfloor,m\beta-\lfloor m\beta\rfloor)_{m\in I^{>0}}$ lie
on, and only on, those portion of the lines $ax+by=t$, where $t$
is any integer, lying within the unit square. Furthermore these
points are dense on these segments.
\end{theorem}

So we have the following corollaries by methods similar to which
represented in \cite[Section 3.5]{N}:

\begin{corollary}({\bf Fact C}) If $I\models B\acute{e}z$ is an $\IP$ for $F$
and suppose $\alpha,\beta$ are positive irrationals such that
$a\alpha^{-1}+b\beta^{-1}=c$ for some $a,b,c\in I$ with $ab<0$ and
$c\neq 0$. Then $N_{\alpha}\cap N_{\beta}$ is a cofinal subset of
$I^{\geqslant 0}$.
\end{corollary}

\begin{corollary}({\bf Fact D}) Suppose $I$ is a B\'ezout {\bf EDR}
and it is an $\IP$ for $F$. Let $\alpha,\beta$ be positive
irrationals such that $a\alpha^{-1}+b\beta^{-1}=c$ for some
$a,b,c\in I^{>0}$ with $c>1$ and $(a,b,c)=1$. Then $N_{\alpha}\cap
N_{\beta}$ is a cofinal subset of $I^{\geqslant 0}$.
\end{corollary}

\noindent We showed that if $a\alpha^{-1}+b\beta^{-1}=1$, then
$N_{\alpha}\cap N_{\beta}={0}$. So if $I\models B\acute{e}z$, the
reminder case is

\begin{center} $\{1,\alpha^{-1},\beta^{-1}\}$ are
linear independent over $Frac(I)$.
\end{center}

\noindent This case is {\it Kronecker's Theorem}. We don't know
whether $\mathbb{P}\vdash Kronecker's\ Th.$ or not. If not, we
must have some $\langle F,I\rangle\models \mathbb{P}+(\neg
Kronecker's\ Th.)$.

\section{Dirichlet's Theorem and Weak Fragments of Arithmetic}
In this section, we prove the Dirichlet's Theorem and consequently
the $\DMO$ property for a nontrivial structure $\langle
F,I\rangle$. Classic proof of Dirichlet's Theorem is based on
${\bf PHP}$. Using this fact, P. D'Aquino proved a weak version of
this theorem, \cite{D}.

\subsection{Weak PHP and Dirichlet's Theorem}
P. D'Aquino studied the theory of Pell equation in $I\Delta_0$.
She used a weak version of ${\bf PHP}$ which is called
$\mathbf{\Delta_0-WPHP}$:
$$ for\ all\ x\ there\ is\ no\ 1-1\ \Delta_0-function\ f\ such\
that\ f:2x\longrightarrow x.$$

\noindent The principle $\mathbf{\Delta_0-WPHP}$ is available in
the theory $I\Delta_0+\Omega_1$, where $\Omega_1$ is
$$\forall x \exists y (x^{\lfloor \log_2{x}\rfloor}=y).$$
The system $I\Delta_0+\Omega_1$ has been widely studied. We know
that
$$IE_1\subset IE_2\subset\cdots I\Delta_0\varsubsetneq I\Delta_0+\Omega_1.$$
P. D'Aquino proved the following version of Dirichlet's Theorem:

\begin{theorem}(\cite[Theorem 3.1]{D}) Let
$\mathcal{M}\vDash I\Delta_0+\Omega_1$, $d\in\mathcal{M}$, $d$ not
a square, $Q>1$, then there are $p,q\in\mathcal{M}$ such that
$|p-\sqrt{d}q|<\frac{1}{Q}$, and $q<2Q$.
\end{theorem}

We will prove a more strong version of Dirichlet's Theorem without
{\it using} ${\bf PHP}$ or any weak version of it in the $IE_1$
system.
\subsection{Farey Series And $IE_{1}$}

First, we define {\it Farey series}. Then we prove some property
of these series. Basic definitions and notations of this
subsection are based on \cite{HW}.

\begin{definition} {\em Suppose $I\vDash GCD$. For an arbitrary
$N\in I^{>0}$, we can define {\it Farey series} $\mathfrak{F}_{N}$
of order $N$ as follows. The Farey series $\mathfrak{F}_{N}$ is
the ascending series of irreducible fractions between $0$ and $1$
whose denominators do not exceed $N$. Thus $\frac{h}{k}\in
\mathfrak{F}_{N}$ if $0\leqslant h\leqslant k\leqslant N,\
(h,k)=1$.}
\end{definition}

We usually suppose $0\in\mathfrak{F}_{N}$. Now, we prove some
important properties of $\mathfrak{F}_{N}$.

\begin{theorem}\label{FT} Suppose $I\models GCD$. If $N\in
I^{>0}$ and $0,1\neq\frac{a}{b}\in \mathfrak{F}_{N}$ and exist
$x_{0},y_{0}\in I$ such that $bx_{0}-ay_{0}=1$, then
\begin{enumerate}
\item there exists unique
successor for $\frac{a}{b}$ in $\mathfrak{F}_{N}$.

(i.e. there exists $\frac{c}{d}\in \mathfrak{F}_{N}$ such that
$\frac{a}{b}<\frac{c}{d}$ and for all $\frac{a}{b}<\frac{m}{n}\in
\mathfrak{F}_{N}$, $\frac{c}{d}\leqslant \frac{m}{n}$.)
\item there exists unique
pre-successor for $\frac{a}{b}$ in $\mathfrak{F}_{N}$.

(i.e. there exists $\frac{c}{d}\in \mathfrak{F}_{N}$ such that
$\frac{a}{b}>\frac{c}{d}$ and for all $\frac{a}{b}>\frac{m}{n}\in
\mathfrak{F}_{N}$, $\frac{c}{d}\geqslant \frac{m}{n}$.)
\end{enumerate}
\end{theorem}

\begin{proof} {\bf {1)}} Since $(x_{0},y_{0})$ is a solution of
$bx-ay=1$, for each $r\in I$, $(x_{0}+ra,y_{0}+rb)$ is also a
solution for $bx-ay=1$. Choose $r$ such that
$N-b<y_{0}+rb\leqslant N$, (we can do it by choose
$r=\lfloor\frac{N-y_{0}}{k}\rfloor$). Now define $x=x_{0}+rb$,
$y=y_{0}+ra$. Therefore, $N-b<y\leqslant N$, $bx=1+ay$. Thus
$x=\frac{a}{b}y+\frac{1}{b}<y+\frac{1}{b}$. Then
$(x,y)=1,x\leqslant y$ and we have
$\frac{x}{y}\in\mathfrak{F}_{N}$.

Note that $\frac{x}{y}=\frac{a}{b}+\frac{1}{ky}>\frac{a}{b}$.
Consequently, $\frac{x}{y}$ appears after $\frac{a}{b}$ in
$\mathfrak{F}_{N}$. If it is not successor of $\frac{a}{b}$, there
exists some $\frac{h}{k}$ between $\frac{a}{b}$ and $\frac{x}{y}$.
So we have \begin{center}
$\frac{x}{y}-\frac{h}{k}=\frac{kx-hy}{ky}\geqslant\frac{1}{ky}$

$\frac{h}{k}-\frac{a}{b}=\frac{bh-ak}{bk}\geqslant\frac{1}{bk}$.
\end{center}
On the other hand, we have

\[\frac{1}{by}=\frac{x}{y}-\frac{a}{b}=(\frac{x}{y}-\frac{h}{k})+(\frac{h}{k}-\frac{a}{b})
\leqslant \frac{1}{ky}+\frac{1}{bk}=\frac{y+b}{bky}.\]

But $y+b>N$. Thus
$\frac{x}{y}-\frac{a}{b}>\frac{N}{bky}\geqslant\frac{1}{by}$. It
is a contradiction.

{\bf {2)}} By the similarly method, we have
$a(-y_{0})-b(-x_{0})=1$. For all $r\in I$,
$y=-y_{0}+rb,x=-x_{0}+ra$ is also a solution for $ay-bx=1$. Now
choose $r=\lfloor\frac{N+y_{0}}{b}\rfloor$. Then
$\frac{x}{y}\in\mathfrak{F}_{N}$. We have
$\frac{x}{y}<\frac{a}{b}$ and moreover
$\frac{a}{b}-\frac{x}{y}=\frac{1}{by}$. By the similar
inequalities, one can prove $\frac{x}{y}$ is a pre-successor of
$\frac{a}{b}$.
\end{proof}

\noindent It seems that the assumptions of Theorem \ref{FT} is
essential, i.e. we have the following claim:

\noindent{\bf Claim} Let $I$ be a $GCD$ domain and $N\in I^{>0}$.
If there exists some $\frac{a}{b}\in \mathfrak{F}_{N}$ such that
for all $x,y\in I$, $ax-by\neq 1$, then $\frac{a}{b}$ has no
successor and pre-successor.

If $I\vDash B\acute{e}z$, then the assumptions of Theorem \ref{FT}
are hold. So every element of $\mathfrak{F}_{N}$ which is not
$0,1$ has successor and pre-successor. For $0$, we have the
successor $\frac{1}{N}$ and for $1$, we have the pre-successor
$\frac{N-1}{N}$. We could prove some properties of
$\mathfrak{F}_{N}$ for these integer parts:

\begin{lemma}\label{LFT} If $I\models B\acute{e}z$, and $N\in
I^{>0}$, then \begin{enumerate}
\item If $\frac{h}{k}$ and $\frac{h'}{k'}$ are two successive
elements of $\mathfrak{F}_{N}$, then $k+k'>N$.
\item No two successive elements of $\mathfrak{F}_{N}$ has the same
denominator.
\item If $\frac{h}{k}$ and $\frac{h'}{k'}$ are two successive
elements of $\mathfrak{F}_{N}$, then $kh'-hk'=1$.
\end{enumerate}
\end{lemma}
\begin{proof} {\bf 1)} The mediant $\frac{h+h'}{k+k'}$ of
$\frac{h}{k}$ and $\frac{h'}{k'}$, falls in the interval
$(\frac{h}{k},\frac{h'}{k'})$. So, if $k+k'\leqslant N$, then
$\frac{h+h'}{k+k'}$ or the reduced format of  it is in
$\mathfrak{F}_{N}$ and it is between $\frac{h}{k}$ and
$\frac{h'}{k'}$.

{\bf 2)} If $k>1$, and $\frac{h'}{k}$ succeeds $\frac{h}{k}$ in
$\mathfrak{F}_{N}$, then $h+1\leqslant h'<k$. But we have
$hk<(h+1)(k-1)$, and therefore $\frac{h}{k-1}<\frac{h+1}{k}$. Then
$\frac{h}{k}<\frac{h}{k-1}<\frac{h+1}{k}\leqslant \frac{h'}{k}$.
But $\frac{h}{k-1}$ comes between $\frac{h}{k}$ and $\frac{h'}{k}$
in $\mathfrak{F}_{N}$, a contradiction.

{\bf 3)} Since $(h,k)=1$, the equation $kx-hy=1$ is soluble in
$I$. If $(x_{0},y_{0})$ is a solution, then $(x_{0}+rh,y_{0}+rk)$
is also a solution for any $r\in I$. We can choose $r$ so that
$N-k<y_{0}+rk\leqslant N$. For this $N-k-y_{0}<rk\leqslant
N-y_{0}$. Then $\frac{N-y_{0}}{k}-1<r\leqslant\frac{N-y_{0}}{k}$.
So $r=\lfloor\frac{N-y_{0}}{k}\rfloor$. Therefore, there is a
solution $(x,y)$ of the equation $kx-hy=1$ such that $(x,y)=1$ and
\[ 0\leqslant N-k<y\leqslant N.\]
Note that $\frac{x}{y}$ is in its lowest terms and $y\leqslant N$
and we have $x=y\frac{h}{k}+\frac{1}{k}<y+\frac{1}{k}$. So
$0<x\leqslant y$. Thus $\frac{x}{y}\in\mathfrak{F}_{N}$. Also, we
have
$\frac{x}{y}=\frac{hy+1}{ky}=\frac{h}{k}+\frac{1}{ky}>\frac{h}{k}$.
So that $\frac{x}{y}$ comes later in $\mathfrak{F}_{N}$ than
$\frac{h}{k}$. If it is not $\frac{h'}{k'}$, it cames later than
$\frac{h'}{k'}$ and
$\frac{x}{y}-\frac{h'}{k'}=\frac{k'x-h'y}{yk'}$. Thus we have
$\frac{x}{y}-\frac{h'}{k'}\geqslant\frac{1}{yk'}$.

While $\frac{h'}{k'}-\frac{h}{k}=\frac{kh'-hk'}{kk'}$, then
$\frac{h'}{k'}-\frac{h}{k}\geqslant\frac{1}{kk'}$. Hence
$\frac{1}{ky}=\frac{kx-hy}{ky}$, which equals to
$\frac{x}{y}-\frac{h}{k}$. But it is less than or equals to
$\frac{1}{yk'}+\frac{1}{kk'}$. The latter is equal to
$\frac{k+y}{kk'y}$. But we have $N-k\leqslant y$, therefore
$y+k>N$. Then
$\frac{k+y}{kk'y}>\frac{N}{kk'y}\geqslant\frac{1}{ky}$. This is a
contradiction and therefore $\frac{x}{y}$ must be $\frac{h'}{k'}$
and $kh'-hk'=1$.
\end{proof}

Suppose that $I\models GCD$. For each $N\in I^{>0}$, we define a
function $\varphi_{N}:\mathfrak{F}_{N}\rightarrow I^{>0}$ by
$\varphi_{N}(\frac{a}{b})=\lfloor N^{2}\frac{a}{b}\rfloor$.
$\varphi_{N}$ is an embedding because if $\frac{a}{b}<\frac{c}{d}$
in $\mathfrak{F}_{N}$, then
$\frac{c}{d}-\frac{a}{b}=\frac{bc-ad}{bd}$. Since $bc-ad\in
I^{>0}$, we have $\frac{bc-ad}{bd}\geqslant\frac{1}{bd}$.
Moreover, $b,d\leqslant N$, then $\frac{c}{d}-\frac{a}{b}
\geqslant\frac{1}{N^{2}}$. So $ N^{2}\frac{a}{b}+1\leqslant
N^{2}\frac{c}{d}$. Then $\varphi_{N}(\frac{a}{b})=\lfloor
N^{2}\frac{a}{b}\rfloor<\lfloor
N^{2}\frac{c}{d}\rfloor=\varphi_{N}(\frac{c}{d})$.

\begin{lemma}\label{FT1} Suppose that $I\vDash IE_{1}$. For each
$m\leqslant N^{2}$, there exists a greatest element of
$\mathfrak{F}_{N}$ such as $\frac{x}{y}$, for which
$\varphi_{N}(\frac{x}{y})<m$.
\end{lemma}
\begin{proof} In fact, we show that $\{n |\ n\in I^{\geqslant
0},\ n\leqslant m,\exists
\frac{x}{y}\in\mathfrak{F}_{N},\varphi_{N}(\frac{x}{y})=n\}$ is a
nonempty bounded $E_{1}$-definable set. Note that
$n=\varphi(\frac{x}{y})$ iff $ny\leqslant N^{2}x<(n+1)y$.
Therefore $x<(n+1)\frac{y}{N^{2}}\leqslant \frac{n+1}{N}$. Thus we
have $x<n$, as a weak inequality. On the other hand, $(x,y)=1$ is
an $E_{1}$-definable sentence.

\noindent Now we define the following $E_{1}$-definable bounded
subset with parameters $N,m$:
\[ \exists x\leqslant mN,\ \exists 0<y\leqslant N,\ x<y, (x,y)=1\
\wedge\ ny\leqslant N^{2}x<(n+1)y\ \wedge n\leqslant m. \]

\noindent The above bounded subset is nonempty, so it has a
greatest element such as $n_{0}$, ( see \cite[(lemma 1.5)]{Sm}).
We have a unique element as $\frac{x}{y}\in\mathfrak{F}_{N}$ with
respect to $n_{0}$. For this element, we have
$\varphi_{N}(\frac{x}{y})=n_{0}\leqslant m$ and $\frac{x}{y}$ is a
greatest element of $\mathfrak{F}_{N}$ for which this property
holds.
\end{proof}

\begin{theorem}{\bf{(Dirichlet's Approximation Lemma in
$IE_{1}$)}} \label{DAL}

Suppose that $I\models IE_{1}$. If $\alpha\in F\setminus Frac(I)$,
for every $Q\in I^{>0}$, there exist $p,q\in I$ with $1\leqslant
q\leqslant Q$, $(p,q)=1$ and
$|\frac{p}{q}-\alpha|\leqslant\frac{1}{qQ}$.
\end{theorem}
\begin{proof} It suffices to prove the result when $0<\alpha<1$.
By Theorem \ref{FT1}, for $m=\lfloor Q^{2}\alpha\rfloor\leqslant
Q^{2}$ there exists a greatest element of $\mathfrak{F}_{Q}$ such
as $\frac{x}{y}$, for which $\varphi_{Q}(\frac{x}{y})<m$. On the
other hand $\frac{m}{Q^{2}}<\alpha<\frac{m+1}{Q^{2}}$ and there
exists at most one element of $\mathfrak{F}_Q$ between
$\frac{m}{Q^{2}}$ and $\frac{m+1}{Q^{2}}$. So the number $\alpha$
lies between two terms of the Farey series $\mathfrak{F}_{Q}$, say
$\frac{p_{1}}{q_{1}}<\alpha<\frac{p_{2}}{q_{2}}$, ( Note that
$\frac{p_{1}}{q_{1}}=\frac{x}{y}$ obtained from Lemma \ref{FT1}
and $\frac{p_{2}}{q_{2}}$ obtained from part (1) of Theorem
\ref{FT} as successor of $\frac{p_{1}}{q_{1}}$). Consider the
mediant $\frac{p_{1}+p_{2}}{q_{1}+q_{2}}$; because this lies
between $\frac{p_{1}}{q_{1}}$ and $\frac{p_{2}}{q_{2}}$ and does
not appear in $\mathfrak{F}_{Q}$, we must have
$q_{1}+q_{2}\geqslant Q+1$. Now $\alpha$ lies in one and only one
of the intervals
$(\frac{p_{1}}{q_{1}},\frac{p_{1}+p_{2}}{q_{1}+q_{2}})$,
$(\frac{p_{1}+p_{2}}{q_{1}+q_{2}},\frac{p_{2}}{q_{2}})$.

If it lies in the first then,
$|\alpha-\frac{p_{1}}{q_{1}}|\leqslant\frac{p_{1}+p_{2}}{q_{1}+q_{2}}-\frac{p_{1}}{q_{1}}$.
The latter is equal to
$\frac{p_{2}q_{1}-q_{2}p_{1}}{q_{1}(q_{1}+q_{2})}$. Since
$p_{2}q_{1}-q_{2}p_{1}=1$,
$\frac{p_{2}q_{1}-q_{2}p_{1}}{q_{1}(q_{1}+q_{2})}=\frac{1}{q_{1}(q_{1}+q_{2})}$.
But we have $q_{1}+q_{2}\geqslant Q+1$, so
$|\alpha-\frac{p_{1}}{q_{1}}|\leqslant\frac{1}{q_{1}(Q+1)}$.
Finally it is less than $\frac{1}{q_{1}Q}$, if we put $p=p_{1}$,
$q=q_{1}$.

Similarly, if it lies in the second, then
$|\alpha-\frac{p_{2}}{q_{2}}|\leqslant
\frac{p_{2}}{q_{2}}-\frac{p_{1}+p_{2}}{q_{1}+q_{2}}$. The latter
is equal to $\frac{p_{2}q_{1}-q_{2}p_{1}}{q_{2}(q_{1}+q_{2})}$.
Since $p_{2}q_{1}-q_{2}p_{1}=1$,
$\frac{p_{2}q_{1}-q_{2}p_{1}}{q_{2}(q_{1}+q_{2})}=
\frac{1}{q_{2}(q_{1}+q_{2})}$. But we have $q_{1}+q_{2}\geqslant
Q+1$, so
$|\alpha-\frac{p_{2}}{q_{2}}|\leqslant\frac{1}{q_{2}(Q+1)}$.
Finally it is less than $\frac{1}{q_{2}Q}$ and we may take
$p=p_{2}$, $q=q_{2}$.
\end{proof}

The above format of Dirichlet's Theorem has some difference by
Dirichlet's Theorem mentioned in Introduction. But we show that
they are the same.
\begin{corollary}\label{CD} Suppose that $I\models IE_{1}$. If $\alpha\in F\setminus Frac(I)$,
for every $Q\in I^{>0}$, there exist some $h,k\in I$ such that
$k>Q$ and
\[|\alpha-\frac{h}{k}|<\frac{1}{k^2}\ .\]
\end{corollary}
\begin{proof} Suppose $0<\alpha<1$. First note that if $N_1<N_2$, the fractions in the
$\mathfrak{F}_{N_2}\setminus\mathfrak{F}_{N_1}$ have denominators
larger than $N_1$.

Fix $Q\in I^{>0}$, and let $\frac{h}{k}\in\mathfrak{F}_Q$ be such
that $|\alpha-\frac{h}{k}|<\frac{1}{kQ}$ and $k\leqslant Q$,
certainly. Thus if $\frac{h}{k}$ and $\frac{h'}{k'}$ are two
successive elements of $\mathfrak{F}_Q$ such that $\alpha$ lies
between them, then set
$\epsilon=min\{|\alpha-\frac{h}{k}|,|\alpha-\frac{h'}{k'}|\}$ and
$N=\lfloor\frac{2}{\epsilon}\rfloor+1$. It is obvious that $N>Q$.
Now consider the set $\mathfrak{F}_N$. So by the first paragraph
of proof, there exist some fractions of $\mathfrak{F}_N$ which lye
certainly between $\alpha$ and $\frac{h}{k}$ and there exist some
fractions of $\mathfrak{F}_N$ which lye certainly between $\alpha$
and $\frac{h'}{k'}$. These fractions have denominators greater
than $Q$. Suppose $\frac{p}{q}$ and $\frac{p'}{q'}$ in
$\mathfrak{F}_N$ such that $\alpha$ lies between them. Then the
required inequality holds with $\frac{h}{k}$ replaced by at least
one $\frac{p}{q}$, $\frac{p+p'}{q+q'}$ and $\frac{p'}{q'}$.
\end{proof}

We can see Dirichlet's Approximation Lemma proves $\mathbb{P}$
property in $IE_{1}$-models, (see proposition \ref{P1}(ii)). So
the structures mentioned in theorem \ref{FT1} are separable.
Moreover Wilmers showed that $IE_1\models B\acute{e}z$, \cite{Wi}.
So we have
\begin{theorem} Suppose that $I\models IE_{1}$. Then
\begin{enumerate}
\item ({\bf Fact C}) If $\alpha,\beta$ be positive irrationals such that
$a\alpha^{-1}+b\beta^{-1}=c$ for some $a,b,c\in I$ with $ab<0$ and
$c\neq 0$. Then  $N_{\alpha}\cap N_{\beta}$ is a cofinal subset of
$I^{\geqslant 0}$.

\item ({\bf Fact D}) If $\alpha,\beta$ be positive irrationals such that
$a\alpha^{-1}+b\beta^{-1}=c$ for some $a,b,c\in I^{>0}$ with $c>1$
and $(a,b,c)=1$. Then $N_{\alpha}\cap N_{\beta}$ is a cofinal
subset of $I^{\geqslant 0}$.

\item ({\bf Fact F'}) Suppose $\rho$ and $\sigma>1$ are
rational. If $N_{\rho}\subseteq N_{\sigma}$, then there exist
$a,b\in I^{\geqslant 0}$ such that
$a\rho^{-1}+b(1-\sigma^{-1})=1$.
\end{enumerate}
\end{theorem}

\noindent Corollary \ref{CD} provides a symmetric rational
approximation for every irrational element $\alpha$:
\[-\frac{1}{q^2}<\alpha-\frac{p}{q}<\frac{1}{q^2}.\]

\noindent In \cite{Se}, B. Segre proved an asymmetric version of
Dirichlet's Theorem. Niven presented a proof using Farey series,
(see \cite[Section 1.3]{N}). In the rest of this section, we will
show that this asymmetric Diophantine approximations Theorem holds
for structures mentioned in theorem \ref{DAL}.

Applying proofs similar to the proof of Corollary \ref{CD}, we
conclude that if $r\in I$ is a positive element, then for all
sufficiently large number $n$, the two fractions $\frac{a}{b}$ and
$\frac{c}{d}$ adjacent to $\alpha$ in $\mathfrak{F}_n$ have
denominators larger than $r$, that is, $b>r$ and $d>r$. We begin
with a preliminary result.

\begin{lemma}
Let $\alpha$ be an irrational and $\tau>0$. If $\frac{a}{b}$ and
$\frac{c}{d}$ are rational numbers with positive denominators such
that $bc-ad=1$ and
\[\frac{a}{b}<\alpha<\frac{c}{d}\ .\]

\noindent Then the following inequalities holds with $\frac{h}{k}$
replaced by at least one of $\frac{a}{b}$, $\frac{a+c}{b+d}$ and
$\frac{c}{d}$:
\[-\frac{1}{\sqrt{(1+4\tau)}k^2}<\alpha-\frac{h}{k}<\frac{\tau}{\sqrt{(1+4\tau)}k^2}\
.\]
\end{lemma}
\begin{proof}
The proof is similar to the proof Lemma 1.8 in \cite{N}.
\end{proof}

Now we can provide an asymmetric version of Dirichlet's Theorem.
\begin{theorem}
Suppose that $I\models IE_{1}$. If $\alpha\in F\setminus Frac(I)$
and $\tau$ is an arbitrary positive element. For each element
$Q\in I^{>0}$, there exists $h,k\in I$ such that $K>Q$ and
\[-\frac{1}{\sqrt{(1+4\tau)}k^2}<\alpha-\frac{h}{k}<\frac{\tau}{\sqrt{(1+4\tau)}k^2}\
.\]
\end{theorem}
The proof is similar to the proof of \cite{N}, Theorem 1.7. So we
only give one interesting corollary of this theorem.
\begin{corollary}{\bf{(Hurwitz's Approximation Lemma in
$IE_{1}$)}}

Suppose that $I\models IE_{1}$. If $\alpha\in F\setminus Frac(I)$,
for every $Q\in I^{>0}$, there exists $\frac{p}{q}$ with $q>Q$,
$(p,q)=1$ and $|\frac{p}{q}-\alpha|<\frac{1}{\sqrt{5}q^2}$.
\end{corollary}
\begin{proof}
It suffices to let $\tau=1$ in the previous theorem.
\end{proof}
\begin{corollary}
Suppose that $I\models IE_{1}$. If $\alpha\in F\setminus Frac(I)$,
for every $Q\in I^{>0}$, there exist $\frac{p_i}{q_i}$ with
$q_i>Q$, $(p_i,q_i)=1$ for $i=1,2$ such that
$0<\frac{p_1}{q_1}-\alpha<\frac{1}{q_1^2}$ and
$0<\alpha-\frac{p_2}{q_2}<\frac{1}{q_2^2}$.
\end{corollary}
\begin{proof}
It suffices to let $\tau=0$ for $\frac{p_1}{q_1}$ and replace
$\alpha$ by $-\alpha$ for $\frac{P_2}{q_2}$.
\end{proof}
\section{Concluding Remarks and Questions}

In this section, we mention some related questions and partial
results.
\subsection{Separable Fields}
By the remark after Lemma \ref{Linf}, we proved the $\mathbb{S}$
property for a wide class of elements of an arbitrary structure
$\langle F,I\rangle$. The remaining case is when

``$\beta$ is irrational and $\rho$ is ratioanl such that
$1<\beta<\rho<2$ and $(m+1)\frac{\rho}{\rho-1}\in I$ for
$m=\lfloor\frac{(\beta-1)(\rho-1)}{\rho-\beta}\rfloor$."

\noindent In this case, $(m+1)\frac{\rho}{\rho-1}\in N_\rho\cap
N_\beta\cap N_{\frac{\rho}{\rho-1}}$. Therefore
$k=\frac{m+1}{\rho-1}\in I$ and $k\rho=(m+1) \frac{\rho}{\rho-1}$.
We have the following claim:

{\bf Claim}- In the above case,
$\lfloor(k+1)\rho\rfloor=k\rho+1\in N_\rho\setminus N_\beta$.

\subsection{Kronecker's Theorem and Farey series}
In the classical case, all implications of Theorem
\ref{irrational} are reversible (see \cite{N}). Nevertheless, one
can show that in the general $\langle F,I\rangle$ context, if the
condition $N_\alpha\cap N_\beta=\{0\}$ implies the existence of
$a,b\in I^{>0}$ with $a\alpha^{-1}+b\beta^{-1}=1$, then all of the
aforementioned implications are reversible. Furthermore, in this
casde, the $\DMO$ property hold. These results depend on
Kronecker's two dimensional $\DMO$ Theorem as appeared in
\cite{N}. It seems that the one dimensional $\DMO$ does not imply
the two dimensional case. It is very interesting to prove
Kronecker's Theorem without the assumption of PHP and only by
using Farey series.

\begin{question} Does there exist any countable model $\langle
F,+,\cdot,<,I\rangle$ satisfying $\DMO$ in which Kronecker's
Theorem fails?
\end{question}

We showed that the $\DMO$ property and the $\mathbb{P}$-condition
are equivalent. Note that the $\DMO$ property is a first order
sentence for $\langle F,+,\cdot,<,I\rangle$. So by the downward
L\"{o}wenheim-Skolem theorem, it suffices to find out the answer to
the following

\begin{question} Does $\DMO$ hold for all countable structures
$\langle F,+,\cdot,<,I\rangle$?
\end{question}

We showed in this paper that the $\mathbb{P}$-condition implies
$\mathbb{S}$.

\begin{question} Can a model $\langle F,+,\cdot,<,I\rangle$
satisfy $\mathbb{S}$ but not the $\mathbb{P}$-condition?
\end{question}

On the other hand, if we can prove the statement of Theorem
\ref{FT1} for $I\vDash B\acute{e}z$, then Dirichlet's
Approximation Lemma will be proved for all $B\acute{e}z$ integer
parts which can be the best result about Dirichlet's Approximation
Lemma.

{\bf Acknowledgment.} This research was in part supported by a
grant from IPM (No. 81030213). It forms part of the first author's
PhD thesis at Tarbiat Modarres University, Tehran, Iran. He would
like to thank Mojtaba Moniri for all his help during the
preparation of the thesis.

\end{document}